\documentclass[12pt]{article}
\usepackage[latin9]{inputenc}
\usepackage{float}
\usepackage{mathtools}
\usepackage{amsmath}
\usepackage{amssymb}
\usepackage{cancel}
\usepackage{stmaryrd}
\usepackage{stackrel}
\usepackage{undertilde}
\usepackage{graphicx}
\usepackage{esint}



\begin{document}
\title{On $p$-adic cascade equations of hydrodynamic type in modeling fully
developed turbulence}
\author{A.\,P.~Zubarev \\
 \textit{ Physics Department, Samara University, } \\
 \textit{ Moskovskoe shosse 34, 443123, Samara, Russia} \\
 \textit{Natural Science Department,} \\
 \textit{Samara State University of Railway Transport,} \\
 \textit{Perviy Bezimyaniy pereulok 18, 443066, Samara, Russia} \\
 e-mail:\:\texttt{apzubarev@mail.ru} }
\maketitle
\begin{abstract}
A $p$-adic hydrodynamic type equation with two integrals of motion
is proposed. It can be considered as a model cascade equation for
energy dissipation in fully developed turbulence. Some of special
cases of the proposed equation are detailed and it is shown that for
a specific choice of parametrs they have stationary solutions that
correspond to the 2/3 Kolmogorov-Obukhov law. Possible further studies
of the proposed model are discussed.

\end{abstract}

\section{Introduction}

$p$-Adic analysis is an adequate tool for describing systems that
have an explicit or hidden hierarchical structure. Over the past 3
decades, the $p$-adic analysis apparatus has been actively used to
model a different class of systems, such as spin glasses, proteins,
models of the genetic code, dynamic systems, socio-economic systems,
etc., (for a review see \cite{ALL}). One of the applications of $p$-adic
analysis is ultrametric modeling of cascade processes of energy dissipation
in fully developed turbulence, outlined in \cite{T1,T2}. The main
idea of work was \cite{T1} to consider the equations of the hierarchical
cascade model in ultrametric space, implemented by the field of $p$-adic
numbers. In this case, such a property of fully developed turbulence
as multifractality is naturally ensured, and it is possible to obtain
the 2/3 Kolmogorov-Obukhov law \cite{K1,K2,F}. A further modification
of the nonlinear ultrametric equation introduced in \cite{T1} was
proposed in work \cite{T2}. It was found that, using the ultrametric
wavelet analysis, one can get a family of exact solutions for this
modified equation.

The equations proposed in \cite{T1,T2} can be considered as some
approximations of the Navier-Stokes equations written in some basis,
after trunking a number of modes and after a possible $p$-adic parameterization.
Note that the question of the existence of integrals of motion of
the proposed equations was not discussed in works \cite{T1,T2}. Nevertheless,
the Navier-Stokes equations in the absence of dissipation and external
forces have two integrals of motion -- energy and hydrodynamic helicity
(or enstrophy in the 2d case). For this reason, $p$-adic equations,
claiming a possible description of cascade processes of energy dissipation
in fully developed turbulence, must have the structure of equations
of hydrodynamic type with two integrals of motion. Recall that the
term "systems of hydrodynamic type" was introduced by Obukhov \cite{Ob_1}
for quadratically non-linear systems that, while keeping only homogeneous
terms, satisfy the following conditions: regularity (conservation
of phase volume) and the existence of at least one quadratic integral
of motion. The first hierarchical models of the hydrodynamic type
with one integral of motion, which describes the cascade process of
energy dissipation and gives the 2/3 law, was proposed by Obukhov
\cite{Ob_1,Ob_2,Ob_3} and Desnyanskii and Novikov \cite{DN}. A hierarchical
model of a hydrodynamic type with two integrals of motion was first
proposed by Gledzer \cite{Gledzer} (see also \cite{DKOC,GDO}), then
developed by Ohkitani and Yamada \cite{OY_1,OY_2} and is now called
the GOY model (see, for example, \cite{KLWB,FDB,BBTT,BJPV}). All
these models belong to the class of the so-called shell models. The
basic idea behind shell models is to divide the spectral space of
thå velocity field into concentric spheres of growing radii $k_{i}=k_{0}\lambda^{i}$,
where $\lambda$ is the parameter characterizing the ratio of the
characteristic scales of turbulent eddies. The set of mods contained
in a single sphere of fixed radius is called the shell. In shell models
only a few mods are retained in each shell. The equations for the
mode components of the velocity field are a system of hydrodynamic
type with two integrals of motion containing only terms that link
only $s$ of the nearest modes. For Obukhov model and Desnyanskii
and Novikov model the value $s=1$, while for GOY model $s=2$ .

In this article, we propose a general $p$-adic model that has the
structure of a hydrodynamic type system with two integrals of motion.
In Section 2, we present the general form of the equations of this
model. In Section 3, we consider some special cases of the general
$p$-adic equation proposed in Section 2 with different values of
$s$ and various forms of the second integral of motion. In this case
the $p$-adic equation is reduced to equations of shell models type.
We show that for certain relations between the parameters the equation
of the model has a stationary solution corresponding to the 2/3 law.

\section{$p$-Adic hydrodynamic type equation with two integrals of motion}

Let us consider the Navier-Stokes equation for an incompressible viscous
fluid, formally written in the form in which the term with pressure
is excluded:
\begin{equation}
\partial_{t}v_{i}=-v_{j}\partial_{j}v_{i}+\left(\partial_{l}\partial_{l}\right)^{-1}\partial_{i}\left(\left(\partial_{j}v_{k}\right)\left(\partial_{k}v_{j}\right)\right)+\nu\partial_{l}\partial_{l}v_{i}+\left(\delta_{ij}-\left(\partial_{l}\partial_{l}\right)^{-1}\partial_{i}\right)f_{j},\;\partial_{j}v_{j}=0.\label{NS}
\end{equation}
In Eq. (\ref{NS}) $i,j,k,l=1,2,3$, $v_{i}=\left(x,t\right)$ is
the velocity field, $\nu$ is the kinematic viscosity, $f_{i}$ is
the density of external forces and summation over repeated indices
is assumed. Let us assume that there exists a certain countable orthonormal
basis of vector functions $\left\{ e_{a}^{(i)}\right\} $, $\intop d^{3}\vec{r}\,e_{a}^{(i)}e_{b}^{(i)}=\delta_{ab}$
($a,b$ are multi-indices) such that $\partial_{i}e_{a}^{(i)}=0$.
To simplify the notation, we consider the functions $\left\{ e_{a}^{(i)}\right\} $to
be real. Then, by decomposing the functions $v_{i}$ and $f_{i}$
into a sum over a set of basis functions $\left\{ e_{a}^{(i)}\right\} $

\[
v_{i}=\sum_{a}v_{a}e_{a}^{(i)},\;f_{i}=\sum_{a}f_{a}e_{a}^{(i)}.
\]
one can write the Eq. (\ref{NS}) in the form of a system of linked
equations

\begin{equation}
\dot{v}_{a}=\sum_{b,c}K_{abc}v_{b}v_{c}-\sum_{j}\nu_{ab}v_{b}+\sum_{j}\phi_{ab}f_{b}.\label{v_g}
\end{equation}
In the absence of viscosity and external forces, the energy per unit
mass $E=\dfrac{1}{2}\intop d^{3}x\,v^{2}$ and the hydrodynamic helicity
$H=\intop d^{3}x\,\varepsilon_{ijk}v_{i}\partial_{j}v_{k}$ are the
integrals of the motion of the Eq. (\ref{NS}). Therefore, the Eqs.
(\ref{v_g}) for $\nu=0$ and $f_{a}=0$ must also have two integrals
of motion, which are

\begin{equation}
E=\dfrac{1}{2}\sum_{a}v_{a}^{2},\label{Energy}
\end{equation}
and
\begin{equation}
H=\dfrac{1}{2}\sum_{a,b}h_{ab}v_{a}v_{b},\;h_{ab}=\intop d^{3}x\,\varepsilon_{ijk}e_{a}^{(i)}\partial_{j}e_{b}^{(k)}.\label{S}
\end{equation}
In this case Eqs. (\ref{v_g}) are a system of hydrodynamic type with
two integrals of motion and its most general form is

\begin{equation}
\dot{v}_{a}=\sum_{b,c,d}\varepsilon_{abd}h_{dc}v_{b}v_{c}-\sum_{b}\nu_{ab}v_{b}+\sum_{b}\phi_{ab}f_{b},\label{v_g_2int}
\end{equation}
where $h_{ab}=h_{ab}$ and $\varepsilon_{abc}$ is totally anti-symmetric
under the exchange of indices.

Next, we represent the $p$-adic parameterization of Eqs. (\ref{v_g_2int}).
First, recall the definition of $p$-adic number. Let $\mathbb{Q}$
be a field of rational numbers and let $p$ be a fixed prime number.
Any rational number $x\ne0$ is uniquely represented as

\begin{equation}
x=\pm p^{n}\frac{a}{b},\label{x_p}
\end{equation}
where $n$ is an integer, and $a$, $b$ are natural numbers that
are not divisible by $p$ and have no common multipliers. The $p$-adic
norm $\left|x\right|_{p}$ of number \textit{$x\in\mathbb{Q}$} is
defined by the equalities $\left|x\right|_{p}=p^{-n}$, $\left|0\right|_{p}=0$.
The field of $p$-adic numbers $\mathbb{Q}_{p}$ is defined as a completion
of the field of rational numbers $\mathbb{Q}$ by $p$-adic norm $\left|x\right|_{p}$.
Any $p$-adic number can be represented as a series converging by
$p$-adic norm:
\[
x=p^{-n}\left(a_{0}+a_{1}p+a_{2}p^{2}+\ldots\right),
\]
The norm on $\mathbb{Q}_{p}$ induces the metric $d(x,y)=\left|x-y\right|_{p}$
which is ultrametric, i.e. satisfies the strong triangle inequality
$\forall\:x,y,z\:d\left(x,y\right)\leq\max\left\{ d\left(y,z\right),d\left(x,z\right)\right\} $.
We will denote: $B_{i}(a)=\{x\in\mathbb{Q}_{p}:\:|x-a|_{p}\leq p^{i}\}$
-- a ball of radius $p^{i}$ centered at point $a$, $S_{i}(a)=\{x\in\mathbb{Q}_{p}:\:|x-a|_{p}=p^{i}\}$
-- a sphere of radius $p^{i}$ centered at point $a$, $B_{i}\equiv B_{i}(0)$,
$S_{i}\equiv S_{i}(a)$, $\mathbb{Z}_{p}\equiv B_{0}$. On $\mathbb{Q}_{p}$
there exists a~unique (up to a factor) Haar measure $d_{p}x$ which
is invariant with respect to translations $d_{p}\left(x+a\right)=d_{p}x$.
We assume that $d_{p}x$ is a full measure; that is,
\begin{equation}
\intop_{\mathbb{Z}_{p}}d_{p}x=1.\label{norm}
\end{equation}
Under this hypothesis the measure $d_{p}x$ is unique. For more information
about $p$-adic numbers, $p$-adic analysis and its applications,
see \cite{VVZ,Sh}.

For what follows, we also define the class $W_{l}^{\alpha}$ $(\alpha\geq0)$
of complex functions $f(x)$ on~$\mathbb{Q}_{p}$ satisfying the
following conditions:

1) $\left|\varphi(x)\right|\le C\left(1+\left|x\right|_{p}^{\alpha}\right)$,
where $C$ is a~real positive number;

2) there exists a~natural number $l$ such that $\varphi\left(x+x'\right)=\varphi\left(x\right)$
for any $x\in\mathbb{Q}_{p}$ and any $x'\in\mathbb{Q}_{p}$, $\left|x'\right|_{p}\le p^{l}$.
A function $\varphi(x)$ satisfying such condition is called locally
constant, and the number $l$ is called the exponent of local constancy
of a function.

For our purposes, we can assume that the number $p$ is a natural
number $p=m>2$. In this case $\mathbb{Q}_{p}$ is a ring of $m$-adic
numbers $\mathbb{Q}_{m}$ with the pseudonorm $\left|x\right|_{m}$,
which also induces on $\mathbb{Q}_{m}$ the ultrametrics $d(x,y)=\left|x-y\right|_{p}$
\cite{DZ}.

Let $v\left(x,t\right)$ be a function on $\mathbb{Q}_{p}\times\mathbb{R}$.
We consider the equation of the following structure
\[
\dfrac{\partial v\left(x,t\right)}{\partial t}=\intop_{\mathbb{Q}_{p}}m\left(x\right)d_{p}y\intop_{\mathbb{Q}_{p}}m\left(z\right)d_{p}z\intop_{\mathbb{Q}_{p}}m\left(z'\right)d_{p}z'\varepsilon\left(x,y,z'\right)h\left(z',z\right)v\left(y,t\right)v\left(z,t\right)
\]
\begin{equation}
-\intop_{\mathbb{Q}_{p}}m\left(y\right)d_{p}y\nu\left(x,y\right)v\left(y,t\right)+f\left(x\right),\label{gen}
\end{equation}
where $m\left(x\right)$ is some function $\mathbb{Q}_{p}\rightarrow\mathbb{R}_{+}$
which is locally integrable with respect to the Haar measure $d_{p}x$,
$\varepsilon\left(x,y,z\right)$ is a completely antisymmetric function
$\mathbb{Q}_{p}\times\mathbb{Q}_{p}\times\mathbb{Q}_{p}\rightarrow\mathbb{R}$
and $h\left(x,y\right)$ is a symmetric function $\mathbb{Q}_{p}\times\mathbb{Q}_{p}\rightarrow\mathbb{R}$.
The Eq. (\ref{gen}) with $\nu\left(x,y\right)=0$, $f\left(x\right)=0$
has the following two integrals of motion:

\begin{equation}
E=\intop_{\mathbb{Q}_{p}}m\left(x\right)d_{p}xv^{2}\left(x,t\right),\:H=\intop_{\mathbb{Q}_{p}}m\left(x\right)d_{p}x\intop_{\mathbb{Q}_{p}}m\left(y\right)d_{p}yh\left(x,y\right)v\left(x,t\right)v\left(y,t\right).\label{E_and_H}
\end{equation}
Eq. (\ref{gen}) is a general $p$-adic equation of hydrodynamic type
with two integrals of motion. We consider special cases of this equation
in the next section.

\section{$p$-Adic cascade equations on $B_{r}$}

The study of the properties of the model based on the Eq. (\ref{gen})
is related to the choice of functions $m\left(x\right)$, $\varepsilon\left(x,y,z\right)$,
$h\left(x,y\right)$, $\lambda\left(x,y\right)$, $f\left(x\right)$.
In this section, we choose the function $m\left(x\right)$), which
determines the integration measure in Eq. (\ref{gen}) as follows
\[
m\left(x\right)=\left(\Omega\left(\left|x\right|_{p}\right)\left(1-p^{-1}\right)+\left(1-\Omega\left(\left|x\right|_{p}\right)\right)\right),
\]
where we use the notation well known in analysis:
\[
\Omega\left(\left|x\right|_{p}p^{-i}\right)=\left\{ \begin{array}{c}
1,\:\left|x\right|_{p}\leq p^{i}\\
0,\:\left|x\right|_{p}>p^{i}
\end{array}\right..
\]
So $v\left(x,t\right)$ is a function on $B_{r}\times\mathbb{R}$
rather than on $\mathbb{Q}_{p}\times\mathbb{R}$. The points $x\in B_{r}\subset\mathbb{Q}_{p}$
parameterize the velocity field modes. We assume that for $y\in\mathbb{Z}_{p}$
the points $x$ and $x+y$ parameterize the same mode of the velocity
field, thus modes are parameterized by $p$-adic balls of unit radius
in $B_{r}$ or, equivalently, by points in $B_{r}/\mathbb{Z}_{p}$.
Also we assume that each mode corresponds to a turbulent eddies and
the scale $l_{i}$ of the eddies corresponding to any unit ball in
$S_{i}$ is determined by the radius of $S_{i}$: $l_{i}=l_{0}p^{-i}$,
where $l_{0}$ is the largest scale of the cascade. In this case,
all modes corresponding to eddies of the same scale $l_{i}$ give
the same contribution to energy and hydrodynamic helicity.

Further, we assume that any eddy of $l_{i}$ can interact only with
eddies of scales $l_{i-s}$, $\ldots$, $l_{i-1}$, $l_{i+1}$,$\ldots$,
$l_{i+s}$. Moreover, the function $\varepsilon\left(x,y,z\right)$
in the Eq. (\ref{gen}), which determines this interaction, has the
form
\[
\varepsilon\left(x,y,z\right)=\sigma\left(x,y\right)\sigma\left(y,z\right)\sigma\left(z,x\right),
\]
\[
\sigma\left(x,y\right)=\sigma_{s}\left(\left|x\right|_{p},\left|y\right|_{p}\right)
\]
\[
=\left|x-y\right|^{-\alpha}\left(\Omega\left(\left|\dfrac{x}{y}\right|_{p}p^{-s}\right)\Omega\left(\left|\dfrac{y}{x}\right|_{p}p\right)-\Omega\left(\left|\dfrac{y}{x}\right|_{p}p^{-s}\right)\Omega\left(\left|\dfrac{x}{y}\right|_{p}p\right)\right),
\]
and $\alpha$ is some real parameter. We also assume that all other
functions in the Eq. (\ref{gen}) are depend only on the $p$-adic
norm: $h\left(x,y\right)=h\left(\left|x\right|_{p},\left|y\right|_{p}\right)$,
$\nu\left(x,y\right)=\nu\left(\left|x\right|_{p},\left|y\right|_{p}\right)$,
$f\left(x\right)=f\left(\left|x\right|_{p}\right)$. Moreover, functions
$f\left(\left|x\right|_{p}\right)\in W_{0}^{0}\left(B_{r}\right)$
and $h\left(\left|x\right|_{p},\left|y\right|_{p}\right),\:\nu\left(\left|x\right|_{p},\left|y\right|_{p}\right)\in W_{0}^{0}\left(B_{r}\right)\times W_{0}^{0}\left(B_{r}\right)$.

Within the framework of the assumptions made, the Eq. (\ref{gen})
has the form
\[
\dfrac{\partial v\left(x,t\right)}{\partial t}=\intop_{\mathbb{Q}_{p}}d\mu(y)\intop_{\mathbb{Q}_{p}}d\mu(z)\intop_{\mathbb{Q}_{p}}d\mu(z')
\]
\[
\times S\left(\left|x\right|_{p},\left|y\right|_{p}\right)S\left(\left|y\right|_{p},\left|z\right|_{p}\right)S\left(\left|z\right|_{p},\left|x\right|_{p}\right)v\left(\left|y\right|_{p},t\right)h\left(\left|z'\right|_{p},\left|z\right|_{p}\right)v\left(\left|z\right|_{p},t\right)
\]
\begin{equation}
-\intop_{\mathbb{Q}_{p}}d\mu(y)\intop_{\mathbb{Q}_{p}}\nu\left(\left|x\right|_{p},\left|y\right|_{p}\right)v\left(y,t\right)+f\left(\left|x\right|_{p}\right).\label{gen_norm}
\end{equation}
where we denote $d\mu(x)\equiv m\left(x\right)d_{p}x$. The Cauchy
problem for the Eq. (\ref{gen_norm}) can be defined in the class
of functions $v\left(x,t\right)=v\left(\left|x\right|_{p},t\right)$
from $W_{0,\mathrm{norm}}^{0}\left(B_{r}\right)\times\mathbb{R}$,
where $W_{0,\mathrm{norm}}^{0}\left(B_{r}\right)$ is the intersection
of the class $W_{0}^{0}\left(B_{r}\right)$ with the class of functions
that depend only on the norm $\left|x\right|_{p}$.

Let $e_{i}\left(x\right)=\Omega\left(\left|x\right|_{p}p^{-i}\right)-\Omega\left(\left|x\right|_{p}p^{-i+1}\right)$
be the characteristic function of the sphere $S_{i}$, $i>0$ and
$e_{0}\left(x\right)=\Omega\left(\left|x\right|_{p}\right)$ be the
characteristic function of $\mathbb{Z}_{p}$. Let us expand the functions
$v\left(x\right)$, $h\left(\left|x\right|_{p}\right)$, $S\left(\left|x\right|_{p},\left|y\right|_{p}\right)$,
$m\left(\left|x\right|_{p}\right)$, $f\left(\left|x\right|_{p}\right)$
in basis $\left\{ e_{0}\left(x\right),e_{1}\left(x\right),\ldots,e_{n}\left(x\right)\right\} $
\[
v\left(x\right)=\sum_{i=0}^{r}V_{i}e_{i}\left(x\right),\:h\left(\left|x\right|_{p},\left|y\right|_{p}\right)=\sum_{i,j=0}^{r}\nu_{i,j}e_{i}\left(x\right)e_{j}\left(y\right),
\]
\[
\nu\left(\left|x\right|_{p},\left|y\right|_{p}\right)=\sum_{i,j=0}^{r}h_{i,j}e_{i}\left(x\right)e_{j}\left(y\right),\:f\left(\left|x\right|_{p}\right)=\sum_{i=0}^{r}f_{i}e_{i}\left(x\right),
\]
\[
S\left(\left|x\right|_{p},\left|y\right|_{p}\right)=\sum_{i=0}^{r}\sum_{j=i-s}^{i-1}p^{-\alpha i}e_{i}\left(x\right)e_{j}\left(y\right)-\sum_{i=0}^{r}\sum_{j=i+1}^{i+s}p^{-\alpha j}e_{i}\left(x\right)e_{j}\left(y\right),
\]
and we accept that $e_{i}\equiv0$ for $i<0$ and $i>r$. Then it
can be shown that the Eq. (\ref{gen_norm}) is equivalent to the following
system of equations for $V_{i}$:
\[
\dot{V}_{i}=h_{0}\left(1-p^{-1}\right)^{3}p^{\left(\gamma+2\right)i}
\]
\[
\times\sum_{j,k,l=0}^{2s}\theta_{s,j}\theta_{s,k}\theta_{s,l}\delta_{j+k+l,3s}p^{-\alpha\left(\max\left\{ s,j\right\} +\max\left\{ s,k\right\} +\max\left\{ s,l\right\} \right)-j-k}V_{i+j-s}\sum_{m=0}^{r}h_{i-k+s,m}p^{m}V_{m}
\]
\begin{equation}
-\sum_{j=0}^{r}\nu_{i,j}V_{j}+f_{i},\label{eq_gen_ODE}
\end{equation}
where
\[
\theta_{i,j}\equiv\left\{ \begin{array}{c}
1,\:i>j\\
-1,\:i<j\\
0,\:i=j
\end{array}\right.
\]
and here and below it is accepted that $V_{i}\equiv0$ for $i<0$
and $i>r$. The Eqs. (\ref{eq_gen_ODE}) for $\nu_{i,j}=0$, $f_{i}=0$
have two integrals of motion:
\begin{equation}
E=\sum_{i=0}^{r}E_{i},\;E_{i}=\left(1-p^{-1}\right)p^{i}V_{i}^{2},\label{sum_E}
\end{equation}
\begin{equation}
H=\left(1-p^{-1}\right)^{2}\sum_{i,j=0}^{r}p^{i}p^{j}h_{i,j}V_{i}V_{j}\label{sum_H}
\end{equation}
Note that the factor $\left(1-p^{-1}\right)p^{i}$ in (\ref{sum_E})
appears because the $p$-adic sphere of radius $p^{i}$ ($i>0)$ is
a union of $\left(1-p^{-1}\right)p^{i}$ $p$-adic balls of unit radius.
Physically, this means that at the hierarchical level $i$, there
is $\left(1-p^{-1}\right)p^{i}$ self-similar eddies, each of which
corresponds to a some $p$-adic unit ball. The velocity component
for each such eddy is equal to $V_{i}$ and therefore the total contribution
of all such eddies to integrals of motion is proportional to the factor
$\left(1-p^{-1}\right)p^{i}V_{i}^{2}$.

Next, we consider a specific choice of function $h\left(\left|x\right|_{p},\left|y\right|_{p}\right)$
of the form
\begin{equation}
h\left(\left|x\right|_{p},\left|y\right|_{p}\right)=h_{0}\left|x\right|_{p}^{\gamma}\Omega\left(\left|\dfrac{x}{y}\right|_{p}\right)\Omega\left(\left|\dfrac{y}{x}\right|_{p}\right),\label{h}
\end{equation}
corresponding to the diagonal form of the second integral of motion
$H$. It is easy to verify that when choosing $s=2$, $\alpha=0$
the Eqs. (\ref{eq_gen_ODE}) go over to the equations
\[
\dot{V}_{i}=h_{0}\left(1-p^{-1}\right)^{3}p^{\left(\gamma+2\right)i}
\]
\[
\times\left(p^{3}\left(p^{\gamma}-p^{2\gamma}\right)V_{i+1}V_{i+2}+\left(p^{\gamma}-p^{-\gamma}\right)V_{i-1}V_{i+1}+p^{-3}\left(p^{-2\gamma}-p^{-\gamma}\right)V_{i-2}V_{i-1}\right)
\]
\begin{equation}
-\sum_{j=0}^{r}\nu_{i,j}V_{j}+f_{i}.\label{Model_s_2}
\end{equation}
Note that the quadratic term in Eqs. (\ref{Model_s_2}) although similar,
do not coincide with corresponding term in the equations of the GOY
model, which have the following general form (see, for example, \cite{KLWB,FDB,BBTT,BJPV}):

\begin{equation}
\dot{v}_{i}=a\lambda^{i}\left(v_{i+1}v_{i+2}-\dfrac{\varepsilon}{\lambda}v_{i-1}v_{i+1}+\dfrac{\varepsilon-1}{\lambda^{2}}v_{i-2}v_{i-1}\right)-\nu_{i}v_{i}+f_{i}.\label{GOY_gen}
\end{equation}
(here $\lambda>1$, $\varepsilon$, $a$ are parameters of model),
since the integrals of motion (\ref{sum_E}) and (\ref{sum_H}) of
Eqs. (\ref{Model_s_2}) are different in form from the integrals of
motion of Eqs. (\ref{GOY_gen}):

\[
E=\sum_{i}V_{i}^{2},\:H=\sum_{i}\left(\varepsilon-1\right)^{-i}V_{i}^{2}.
\]
Note that the correspondence between (\ref{Model_s_2}) and (\ref{GOY_gen})
takes place for $\gamma=-\dfrac{1}{2}$, $p=\lambda$ and after a
scaling $V_{i}=p^{-\frac{1}{2}i}v_{i}.$ One can show that with a
suitable choice of the coefficients $\nu_{i,j}$ and $f_{i}$ the
system of Eqs. (\ref{Model_s_2}) has a stationary solution

\begin{equation}
V_{i}=cp^{-\frac{5}{6}i}\label{st_sol}
\end{equation}
if one of the conditions (\ref{or})
\begin{equation}
\gamma=-\dfrac{1}{2}\:\mathrm{or}\:\gamma=-\dfrac{1}{4}\label{or}
\end{equation}
is satisfied. The stationary solution (\ref{st_sol}) leads to the
following formula for the total energy of the modes of the $i$-th
level

\[
E_{i}\sim p^{i}V_{i}^{2}\sim p^{-\frac{2}{3}i},
\]
which corresponds to the 2/3 law.

We now consider the Eqs. (\ref{eq_gen_ODE}) when choosing (\ref{h}),
$s=3$ and $\alpha=0$. In this case, the Eqs. (\ref{eq_gen_ODE})
have the form
\[
\dot{V}_{i}=h_{0}\left(1-p^{-1}\right)^{3}p^{\left(\gamma+2\right)i}
\]
\[
\times\left(p^{5}\left(p^{3\gamma}-p^{2\gamma}\right)V_{i+3}V_{i+2}+p^{4}\left(p^{3\gamma}-p^{\gamma}\right)V_{i+3}V_{i+1}\right.
\]
\[
+p^{3}\left(p^{2\gamma}-p^{\gamma}\right)V_{i+2}V_{i+1}+p\left(p^{-\gamma}-p^{2\gamma}\right)V_{i-1}V_{i+2}
\]
\[
+\left(p^{-\gamma}-p^{\gamma}\right)V_{i-1}V_{i+1}+p^{-1}\left(p^{-2\gamma}-p^{\gamma}\right)V_{i-2}V_{i+1}
\]
\[
+p^{-3}\left(p^{-\gamma}-p^{-2\gamma}\right)V_{i-1}V_{i-2}+p^{-4}\left(p^{-\gamma}-p^{-3\gamma}\right)V_{i-1}V_{i-3}
\]

\begin{equation}
\left.+p^{-5}\left(p^{-2\gamma}-p^{-3\gamma}\right)V_{i-2}V_{i-3}\right),\label{Model_s_3}
\end{equation}
where the terms corresponding to dissipation and external forces are
omitted. It can be shown that a model with a quadratic term of the
form (\ref{Model_s_3}) for $\gamma=\dfrac{5}{2}$ and $\gamma=\dfrac{5}{4}$
also has a stationary solution (\ref{st_sol}), which corresponds
to the 2/3 law.

In conclusion, we present the case of a model with an off-diagonal
second integral of motion $H$. The following choice
\[
h\left(\left|x\right|_{p},\left|y\right|_{p}\right)=\dfrac{1}{2}h_{0}\left|x\right|^{\gamma}\left(\Omega\left(\left|\dfrac{x}{y}\right|_{p}p^{-1}\right)\Omega\left(\left|\dfrac{y}{x}\right|_{p}p\right)+\Omega\left(\left|\dfrac{y}{x}\right|_{p}p^{-1}\right)\Omega\left(\left|\dfrac{x}{y}\right|_{p}p\right)\right),
\]
leads to the following expression for $H$:
\[
H=h_{0}p^{-1}\sum_{i=0}^{r}p^{\left(\gamma+2\right)i}V_{i}V_{i-1}.
\]
In this case, from (\ref{gen_norm}) we get the equations with qudratic
term of the form
\[
\dot{V}_{i}=h_{0}\left(1-p^{-1}\right)^{3}p^{\left(\gamma+3\right)i}
\]
\[
\times\sum_{j,k,l=0}^{2s}\theta_{s,j}\theta_{s,k}\theta_{s,l}\delta_{j+k+l,3s}p^{-\alpha\left(\max\left\{ s,j\right\} +\max\left\{ s,k\right\} +\max\left\{ s,l\right\} \right)+j-i}
\]
\begin{equation}
\times V_{i+j-s}p^{\left(\gamma+1\right)\left(-k+s\right)}\left(p^{-1}V_{i-k+s-1}+\theta_{i-k+s}p^{\gamma+1}V_{i-k+s+1}\right),\label{dV_nd}
\end{equation}
where
\[
\theta_{i}\equiv\left\{ \begin{array}{c}
1,\:i\geq0\\
0,\:i<j
\end{array}\right..
\]
In the case $s=2$ Eqs. (\ref{dV_nd}) have the form
\[
\dot{V}_{i}=h_{0}\left(1-p^{-1}\right)^{3}p^{\left(\gamma+3\right)i}
\]
\[
\times\left[p^{-3}\left(p^{-\gamma-2}V_{i-2}V_{i-2}+V_{i-2}V_{i}-p^{-2\gamma-3}V_{i-1}V_{i-3}-\left(1-\delta_{i,1}\right)p^{-\gamma-1}V_{i-1}V_{i-1}\right)\right.
\]
\[
+p^{-\gamma-2}V_{i+1}V_{i-2}+\left(1-\delta_{i,0}\right)V_{i+1}V_{i}-\left(1-\delta_{i,n}\right)p^{\gamma}V_{i-1}V_{i}-p^{2\gamma+2}V_{i-1}V_{i+2}
\]
\begin{equation}
+\left.p^{3}\left(\left(1-\delta_{i,n-1}\right)p^{2\gamma+1}V_{i+1}V_{i+1}+p^{3\gamma+3}V_{i+1}V_{i+3}-p^{\gamma}V_{i+2}V_{i}-p^{2\gamma+2}V_{i+2}V_{i+2}\right)\right].\label{model_nd}
\end{equation}
It also can be shown that model (\ref{model_nd}) has a stationary
solution corresponding to the 2/3 law for $\gamma=-1$, $\gamma=-\dfrac{3}{2}$
and $\gamma=-\dfrac{5}{4}$.

\section{Concluding remarks}

The main goal of this short article is to give another important direction
in the application of $p$-adic analysis, namely, the parameterization
of hydrodynamic-type systems with two integrals of motion. Such an
equation was written in the form (\ref{gen}) and in the future it
can be used to study new types of models of cascade processes of energy
dissipation in fully developed turbulence. It was shown that with
a special choice of functions included in this equation, it leads
to a class of shell-type equations of the form (\ref{eq_gen_ODE}).
It was also shown that such a class of equations for a specific choice
of parameters has a stationary solutions that is consistent with Kolmogorov-Obukhov
law. At the same time, we are not discussing problems related to the
stability of stationary solutions of the models under consideration,
and this is one of the direction of future studies. In this connection,
it is also of interest to study different models with off-diagonal
forms of the second integral of motion.

The class of equations considered in Section 3 has discrete scale
invariance, since the norm of the $p$-adic number trivially has the
scaling property: $\left|px\right|_{p}=p^{-1}\left|x\right|_{p}$.
However, these equations (\ref{gen_norm}) are not translationally
invariant on $\mathbb{Q}_{p}$. Note that the general equation (\ref{gen})
allows constructing translationally invariant hierarchical models
that are determined by the choice of the functions $\varepsilon\left(x,y,z\right)$
and $\nu\left(x,y\right)$, incoming in this equation. For example,
the translation-invariant model is obviously obtained by choosing
\[
\varepsilon\left(x,y,z\right)=\sigma\left(x-y\right)\sigma\left(y-z\right)\sigma\left(z-x\right),\;h\left(x,y\right)=h\left(\left|x-y\right|_{p}\right)
\]
where $\sigma\left(x\right)$ is the odd function $\mathbb{Q}_{p}\rightarrow\mathbb{R}$:
$\sigma\left(x\right)=\sigma\left(-x\right)$. Such a class of models
may also be of potential interest for future studies.

\end{document}